\documentclass[12pt,reqno]{amsart}
\usepackage{amsmath}
\usepackage{amsthm}
\usepackage{mathrsfs}
\usepackage{amsfonts} 
\usepackage{amssymb} 
\usepackage{epsfig}

\numberwithin{equation}{section}

\allowdisplaybreaks

\newtheorem{theorem}{Theorem}[section]

\newtheorem{proposition}[theorem]{Proposition}
\newtheorem{lemma}[theorem]{Lemma}

\newtheorem{definition}[theorem]{Definition}

\renewcommand{\epsilon}{\varepsilon}

\newcommand{\R}{\mathbb{R}}
\newcommand{\N}{\mathbb{N}}
\newcommand{\Z}{\mathbb{Z}}
 
\newcommand{\PR}{\mathbb{P}}

\newcommand{\ladder}{\mathcal{B}}

\newcommand{\degg}{d}

\newcommand{\tint}{\mathscr{T}}

\newcommand{\tree}{\mathcal{T}}
\newcommand{\reach}{H}
\newcommand{\added}{\mathcal{A}}
\newcommand{\addedn}{\mathcal{A}_n}

\newcommand{\cross}{\mathsf{C}}
\newcommand{\double}{\mathcal{M}}
\newcommand{\doublephi}{\mathcal{M}_{\phi}}
\newcommand{\single}{\mathcal{S}}
\newcommand{\singlephi}{\mathcal{S}_{\phi}}

\newcommand{\singleedge}{e_{\mathsf{BN}}}
\newcommand{\singlebar}{b_{\mathsf{BN}}}
\newcommand{\singlebarn}{b_{\mathsf{BN}_n}}
\newcommand{\nonesc}{\mathsf{NoEsc}}
\newcommand{\nonescn}{\mathsf{NoEsc}_n}

\newcommand{\pext}{\partial_{\rm ext}}
\newcommand{\scv}{\mathcal{V}}

\newcommand{\seev}{\mathsf{BN}}
\newcommand{\nse}{\seev^c}

\newcommand{\pinf}{p_\infty}

\newcommand{\condladderm}{\PR^{\cross_m \cap \seev_m^c}_{t,(\ladder,\added_m)}}

\newcommand{\ubl}{{\rm UnTouch}}

\newcommand{\forced}{{\rm Found}}

\newcommand{\laddergiven}{\mathcal{B}'}
\newcommand{\viloc}{{\rm ViLoc}}

\newcommand{\onpiv}{\mathsf{P}^+}
\newcommand{\offpiv}{\mathsf{P}^-}

\newcommand{\treen}{\tree_n}

\newcommand{\desctree}[1]{\tree_{[#1]}}
\newcommand{\abovetree}[1]{\tree^{[#1]}}

\newcommand{\hit}{H}

\newcommand{\numericboundstrong}{375}

\newcommand{\nta}{A}

\newcommand{\cardin}[1]{\vert #1 \vert}

\newcommand{\akn}{A_{n,k}}
\newcommand{\aon}{A_{n,0}}

\def\fff#1{&{{\pageref{#1}}}\cr}
\def\hfff#1{\label{#1}}
\title[Sharp phase transition in the random stirring model]{Sharp phase transition in \\ the random stirring model on trees
}
\date{}

\author[A.~Hammond]{Alan HAMMOND}

\address{ Department of Statistics, 
           University of Oxford,
                  1 South Parks Road,
                  Oxford, OX1 3TG, U.K. } 
\email{hammond@stats.ox.ac.uk}

\thanks{Dept. of Statistics, Oxford University. Supported by  U.K. EPSRC grant EP/I004378/1.} 

\begin{document}

\begin{abstract}
We establish that the phase transition for infinite cycles in the random stirring model on an infinite regular tree of high degree is sharp. That is, 
we prove that there exists $d_0$ such that, for any $\degg \geq d_0$,
the set of parameter values at which the random stirring model on the rooted regular tree with offspring degree $\degg$ almost surely contains an infinite cycle consists of a semi-infinite interval. The critical point at the left-hand end of this interval is at least $\degg^{-1} + \tfrac{1}{2}\degg^{-2}$ and at most $\degg^{-1} + 2\degg^{-2}$. 
\end{abstract}

\maketitle

%\tableofcontents

%\newpage 

\vspace{-7mm}

\begin{section}{Introduction}

Suppose given a graph $G = \big( V(G),E(G) \big)$.   To each edge $e \in E(G)$ is associated an independent Poisson process of rate one on $[0,\infty)$. The random stirring model on $G$ is a stochastic process $\sigma$ defined on $[0,\infty)$ and taking values in permutations of $V(G)$. The initial condition $\sigma_0$ is the identity permutation. 
As time $t$ increases, on each occasion that a point $(e,t) \in E(G) \times [0,\infty)$ of one of the Poisson processes is encountered, $\sigma$ is instantaneously modified by composing with the transposition of the two vertices incident to the edge~$e$. 

Let $\degg \geq 2$. Let $\tree$ \hfff{tree} 
denote the rooted regular tree of offspring degree $\degg$, and let $\phi$ \hfff{phi}
denote the root of $\tree$. Our main theorem shows that, if $\degg$ is high, the random stirring model on $\tree$ has a critical value for the transition to infinite cycles.
The theorem confirms 
for such trees a conjecture made in the 1990s by B{\'a}lint T{\'o}th that a wide variety of transitive infinite graphs should exhibit such a critical point; to the author's best knowledge, the conjecture for regular trees  first appeared in print as Conjecture $9$ of \cite{angel}.
\begin{theorem}\label{thm}
Suppose that $\degg \geq 764$. There exists $T_c(\degg) \in (0,\infty)$ such that number of vertices in the cycle of $\phi$ in $\sigma_t$ is almost surely finite if $t < T_c(\degg)$ and is infinite with positive probability if $t > T_c(\degg)$. For such $\degg$, $T_c(\degg) \in \big[ \degg^{-1} + \tfrac{1}{2}\degg^{-2}, \degg^{-1} + 2\degg^{-2} \big]$. 

For each $\epsilon > 0$, there exists $d' \in \N$ such that for $\degg \geq d'$, $T_c(\degg)$ exists and satisfies  $T_c(\degg) \in \big[ \degg^{-1} + \tfrac{1}{2}\degg^{-2}, \degg^{-1} + \big( \tfrac{7}{6} + \epsilon  \big) \degg^{-2}  \big]$.  
\end{theorem}

\subsection{Glossary of notation}
\small{Here we list alphabetically the notation which is commonly used in the article. 
%The list is in Roman alphabetical ordering where, for example, the English spelling of Greek letter names is substituted for those letters. 
A summarizing phrase is provided for each item, as well as the page number at which the concept is introduced.}

%\vspace{-4mm}

\bigskip
\def\qq{&}
%\small{
\begin{center}
\halign{
#\quad\hfill&#\quad\hfill&\quad\hfill#\cr
$\added$ \qq the added bar, with uniform law on $E(\treen) \times [0,1)$ \fff{added}
$\ladder$ \qq bar collection with Poisson-$t$ law on $E(\tree) \times [0,1)$ under $\PR_t$ \fff{ladder}
bar \qq element of $E(\tree) \times [0,1)$ \fff{bar}
$\singlebar$ \qq bottleneck bar, the unique bar in $\ladder$ supported on $\singleedge$ \fff{singlebar}
$\seev$ \qq the event that $\singleedge$ exists \fff{seev}
$b^+,b^-$ \qq the parent and child joints of a bar $b$ \fff{uljoint}
$\cross$ \qq the crossing event: $X^\ladder[0,\reach_n^\ladder] \cap \{ \added^+,\added^- \} \not= \emptyset$  \fff{cross}
%$\lowadded$ \qq close-to-boundary event, $\seev \setminus \highadded$ \fff{lowadded}
$\pext G$ \qq the exterior boundary (a set of edges) of $G \subseteq E(\tree)$ \fff{pext}
$E(b)$ \qq the edge on which the bar $b$ is supported \fff{edgeb}
$\singleedge$ \qq bottleneck edge: last $e \in E(P_{\phi,E(\added)^+})$ supporting unique bar in $\ladder$ \fff{singleedge} 
$\mathcal{E}_i$ \qq $\{ e \in E(\tree): d(\phi,e^+) = i \}$ \fff{edgeset}
$e^+,e^-$ \qq the parent and child endpoint vertices of an edge $e \in E(\tree)$ \fff{ulvertex}
%$\highadded$ \qq far-from-boundary event, $\big\{ d(\phi,\singleedge^-) \leq n - \ncutoff \big\}$  \fff{highadded}
%$\forced_t$ \qq the set of bars in $\ladder$ crossed by $X^\ladder$ during $[0,t]$ \fff{forcedname}
 $\hit_\added^\ladder$ \qq $\inf \big\{ s> 0: X^\ladder(s) \in \{ \added^+,\added^- \} \big\}$\fff{hitadded}
%$\high$ \qq $\big\{ e \in E(\treen): d \big( \phi,e^- \big) \leq n  - \ncutoff \big\}$ \fff{high}
$\reach_n^\ladder$ \qq the hitting time of $\scv_n \times [0,1)$ by $X^\ladder$ \fff{hittimen} 
$\double_v$ \qq the multi-cluster (a set of edges) associated to $v \in V(\tree)$ \fff{multicl}
%$\ncutoff$ \qq parameter for edges to be close to boundary, set in Proposition~\ref{propineqone} \fff{cutoff}
$\nonesc$ \qq the non-escape event: $X^\ladder_{\singlebar^+}$ visits $(\phi,0)$ before $\scv_n \times [0,1)$ \fff{nonesc} 
%$\nobar$ \qq the event that $\ladder \cap \big( \mathcal{E}_0 \times [0,1) \big) = \emptyset$  \fff{nobar}
%$\PR_{t,\ladder}$ \qq the marginal law of $\ladder$ under $\PR_t$, i.e., Poisson-$t$ \fff{prtladder}
$\phi$ \qq the root of $\tree$ \fff{phi}
$p_\infty$ \qq $\PR_t$-probability that $X^\ladder$ never returns to $(\phi,0)$ \fff{pinf}
$\offpiv$ \qq 
the off-pivotal event
$\{ \reach_n^\ladder < \infty \} \cap \{ \reach_n^{\ladder \cup \added} = \infty \}$ 
 \fff{offpiv}
$p_n$ \qq $\PR_t(\reach_n^\ladder < \infty)$ \fff{pn}
pole at $v$ \qq for $v \in V(\tree)$, the set $\{ v \} \times [0,1)$ \fff{pole}
$P_{\phi,v}$ \qq the path in $\tree$ from $\phi$ to $v \in V(\tree)$ \fff{pathv}
$\onpiv$ \qq the on-pivotal event 
$\{ \reach_n^\ladder = \infty \} \cap \{ \reach_n^{\ladder \cup \added} < \infty \}$ \fff{onpiv}
$\single_\phi$ \qq $\{ e \in \pext \double_\phi:  \textrm{$e$ supports a bar in $\ladder$} \}$  \fff{singlephi}
$\tree$ \qq the rooted regular tree with offspring degree $\degg$ \fff{tree}
$\tau$ \qq $t \degg$ \fff{tau}
$\treen$ \qq the subtree of $\tree$ induced by vertices at distance at most $n$ from $\phi$ \fff{treen}
$\desctree{v}$ \qq the descendent tree of $v \in V(\tree)$ 
\fff{desctree}
$\abovetree{v}$ \qq the tree above  $v \in V(\tree)$: $\tree$ after $\desctree{v}$ is excised \fff{abovetree} 
%$\ubl_t$ \qq {\em untouched} bars: $\big\{ b \in E(\tree) \times [0,1): \{ b^+,b^- \} \cap X^\ladder[0,t] = \emptyset \big\}$ \fff{untouched}
$\scv_i$ \qq $\{ v \in V(\tree): d(\phi,v) = i \}$ \fff{scv}
$\viloc_n(\ladder)$ \qq $\{ b \in E(\treen) \times [0,1): \{b^+,b^-\} \cap X^\ladder[0,\reach_n^\ladder] \not= \emptyset, 
E(P_{\phi,E(b)^+}) \subseteq \doublephi \}$  \fff{viloc}
$X^\ladder, X^{\ladder \cup \added}$ \qq shorthand for $X^\ladder_{(\phi,0)}$ and  $X^{\ladder \cup \{ \added \}}$ \fff{meander}
$X^\ladder_{(v,s)}$ \qq cyclic-time random meander from $(v,s) \in V(\tree) \times [0,1)$ \fff{meandergen}
%$^\y\ladder$ \qq the vertex component of $X^\ladder$ \fff{vertexmeander} 
}\end{center}
%}
%\normalsize

%\submit In fact, we will prove the first assertion of Theorem \ref{thm} for $\degg \geq \numericbound$. A few technical refinements yield the conclusion for $\degg \geq \numericboundstrong$. These are furnished in an appendix to the online counterpart \cite{hammondtwo} of this article. 

\newpage

As \cite{angel} mentions, on a regular tree, it is simple to see that, for each $t \in [0,\infty)$, there being positive probability that the cycle of $\phi$ under $\sigma_t$ is infinite is equivalent to the almost sure existence of some infinite cycle under $\sigma_t$. 

\subsection{Cyclic-time random meander and walk}

Our analysis of the random stirring model exploits a closely related dependent random walk which was used in the proof of ~\cite[Theorem 1]{toth} and which was called the cyclic-time random walk  in~\cite{angel}. We now introduce some notation and define this walk.

For each edge $e \in E(\tree)$, the incident vertex of $e$ closer to $\phi$ will be called the parent vertex and will be denoted by $e^+$; the other, called the child vertex and labelled $e^-$. \hfff{ulvertex}
  
For convenience, suppose that $\tree$ is embedded in $\R^2$, so that each element of $V(\tree)$ is identified with a point in $\R^2$ and each element $e \in E(\tree)$ with the line segment $[v_1,v_2] \subseteq \R^2$ where $e = (v_1,v_2)$ for $v_1,v_2 \in V(\tree)$. For each $v \in V(\tree)$, let the pole at $v$,\hfff{pole} $\{ v \} \times [0,1) \subseteq \R^3$, denote the unit line segment that rises vertically from $v$. Elements of $E(\tree) \times [0,1)$ will be called bars.\hfff{bar} A bar $b = (e,h)$ is said to be supported on the edge $e$ and to have height $h$; we also record the edge on which $b$ is supported as $E(b)$.\hfff{edgeb} Note that the bar $(e,h)$ is a horizontal line segment which intersects the poles at $e^+$ and $e^-$; the intersection points $(e^+,h)$ and $(e^-,h)$ will be called the parent and child joints of $(e,h)$. \hfff{uljoint} 

The bar set $E(\tree) \times [0,1)$ carries the product of counting and Lebesgue measure on its components. (As a shorthand, we will refer to this product measure simply as Lebesgue measure.)

Let $\laddergiven \subseteq E(\tree) \times [0,1)$ be a collection of bars which is locally finite in the sense that each $e \in E(\tree)$ supports only finitely many elements of $\laddergiven$. Cyclic-time random meander $X^{\laddergiven}_{(v,h)}:[0,\infty) \to V(\tree) \times [0,1)$,\hfff{meandergen} among $\laddergiven$ and with initial condition $(v,h) \in V(\tree) \times [0,1)$, is the following process. First, $X^\ladder_{(v,h)}(0) = (v,h)$; the process then rises at unit speed on the pole at $v$ until either it reaches $(v,1)$, when it jumps to $(v,0)$, or until it reaches the joint of a bar in $\laddergiven$, when it jumps to the other joint of this bar. After either of these events, $X^{\laddergiven}_{(v,h)}$ continues by iterating the same rule, until it is defined on all of $[0,\infty)$. The process is chosen to be right-continuous with left limits. Note that this choice implies that, if $(v,h)$ is the joint of a bar $b$ in $\laddergiven$, then $X^{\laddergiven}_{(v,h)}$ remains at the pole at $v$ at small times, rather than crossing $b$ at time zero.   
We abbreviate $X^{\laddergiven}$ \hfff{meander} for  $X^{\laddergiven}_{(\phi,0)}$. (There are locally finite choices of $\laddergiven$ for which these rules fail to define $X^{\laddergiven}_{(v,t)}$ on all of $[0,\infty)$. It is a simple matter to verify that this difficulty does not arise in the case that is relevant to us and which we now discuss.)

Let $s \in (0,\infty)$. We will refer to the Poisson law on bar collections of density $s$ with respect to Lebesgue measure on $E(\tree) \times [0,1)$ as the Poisson-$s$ law. Let $\big\{ \ladder_s : s \geq 0 \big\}$
be a coupled collection of random bar collections, where $\ladder_s$ has the Poisson-$s$ law for each $s \in [0,\infty)$ and $\ladder_s \subseteq \ladder_{s'}$ whenever $0 \leq s < s' < \infty$. Define $\sigma_t: V(\tree) \to V(\tree)$ by setting $\sigma_t(v)$ equal to the vertex component of $X^{\ladder_t}_{(v,0)}(1)$, and note that $\sigma_t$ is a random permutation of $V(\tree)$. With this notation, the random stirring process on $\tree$ is the stochastic process, mapping $[0,\infty)$ to permutations of $V(\tree)$, given by $s \to \sigma_s$.

We now fix $t \in (0,\infty)$ and write $\PR_t$ for a probability measure carrying a bar collection $\ladder \subseteq E(\tree) \times [0,1)$ \hfff{ladder} having Poisson-$t$ law. 
Cyclic-time random meander with parameter $t$ is the random process $X^{\ladder}$. 

Cyclic-time random walk (begun at $\phi$) is the vertex-valued process given by projecting $X^\ladder:[0,\infty) \to V(\tree) \times [0,1)$ onto $V(\tree)$. 
%We denote it by $\y^\ladder$.\hfff{vertexmeander} 
%writing $\y^\ladder_{(v,h)}$ for the $V(\tree)$-projection of $X^\ladder_{(v,h)}$. 
(In fact, under our definition, cyclic-time random walk moves at a rate which is a factor of $t$ greater than it does under the definition in \cite{angel}.)  We will discuss cyclic-time random meander rather than walk, and will refer to $X^\ladder$ in shorthand as a meander. See Figure~\ref{cycleexample} for an illustration.  

\subsection{Different perspectives on the random stirring model and other spatial random permutations}

In 1953, Feynman \cite{feynman} wrote the quantum-mechanical partition function for helium as a sum over the energy associated to certain interacting Brownian particles that may interchange their positions over a finite-time interval. He argued that the $\lambda$-transition undergone by the gas at low temperature is reflected by the appearance of large cycles in a measure on permutations naturally associated to this representation of the partition function. 

The random stirring (or random interchange) model was introduced in \cite{harris}. 
In this model also, the conjectured phase transition to lengthy cycles has a physical importance, since it is intimately connected to the off-diagonal long-range order anticipated for the  spin-$1/2$ isotropic quantum Heisenberg ferromagnet at very low temperature:
 B{\'a}lint T{\'o}th in \cite{toth} gave a representation of the partition function for this ferromagnet in terms of the random stirring model. 
 The lecture notes \cite{guw} contain an overview of this topic. The phase transition to infinite cycles proved in Theorem~\ref{thm} is expected to have a counterpart for the Euclidean lattice $\Z^d$ for $d \geq 3$. The author learnt of this question first after it was posed by B{\'a}lint T{\'o}th.

 Recent mathematical progress on the random stirring model 
includes the resolution of Aldous' conjecture identifying its spectral gap \cite{clr}, and a formula for the probability that the random permutation consists of a single cycle \cite{alonkozma}.
%In work addressing what will also be the object of our attention,  
%Omer Angel \cite{angel} has proved that, on a regular tree of degree at least five, and for a certain interval of values of $T$, the random interchange process with parameter $T$ on the tree has infinite cycles almost surely.

The emergence of a giant component under percolation on the complete graph as the percolation parameter increases 
through values near $1/n$ has been intensively studied. Oded Schramm~\cite{comprantrans} showed that this transition is accompanied by the appearance of large-scale cycles in the associated random stirring model:
in the composition of $(1 + \epsilon)n$ independent uniform transposition on a given $n$-set, 
there exists a giant component of edges transposed at least once, of some density $\theta(\epsilon) \in (0,1)$; when the cycle lengths in  this random permutation are normalized by $\theta(\epsilon)n$ and listed in decreasing order, they converge in law to the  Poisson-Dirichlet distribution with parameter one.
Nathana{\"e}l Berestycki \cite{berestycki} has given a short proof that a cycle exists of size $\Theta(n)$ when $(1 + \epsilon)n$ transpositions are made.

\begin{figure}
\centering\epsfig{file=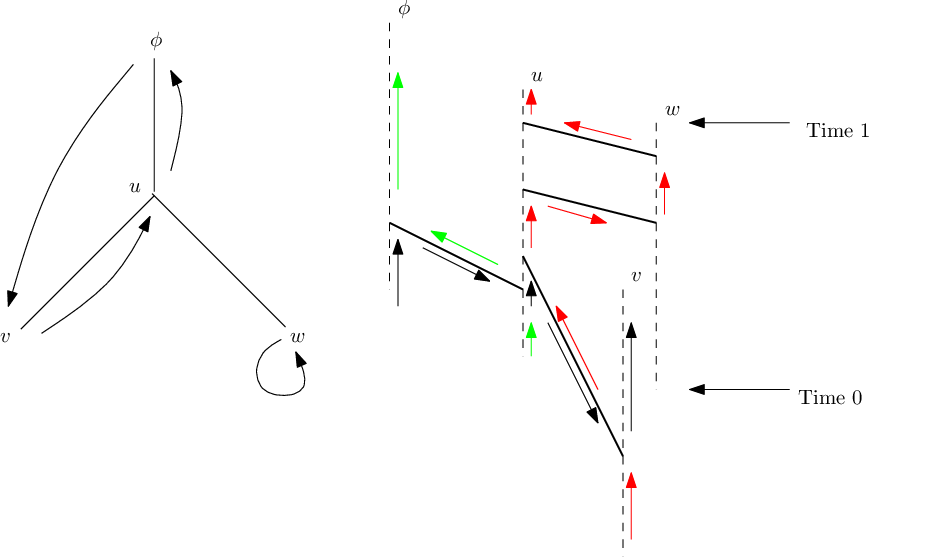, width=14cm}
\caption{For the graph shown on the left, cyclic-time random meander $X^\ladder$ departing from $(\phi,0)$ is illustrated on the right. The right-hand sketch depicts a construction in $\R^3$ in which the poles associated to vertices are the vertical dashed lines and the bars in $\ladder$ are the horizontal black lines.
 Assume that there  are no bars in $\ladder$ supported on edges that  connect vertices $v$ and $w$ of $\phi$ to their offspring. The trajectory of the meander from $(\phi,0)$ is divided into three intervals of unit duration, at the end of which, the meander returns to $(\phi,0)$. These three sub-trajectories are indicated in black, red and green in the right-hand sketch. As the left-hand sketch shows, the cycle of $\phi$ in the associated permutation %(which is $\phi \to v \to u \to \phi$) 
thus has three elements.}\label{cycleexample}
\end{figure}

\subsection{Monotonicity near the transition}

For any given graph $G$ on which the random stirring model is well-defined, 
let $\tint^G$ denote the set of  $t > 0$ such that the random stirring process on $G$ at parameter $t$ has infinite cycles almost surely. 
Note that $t \not\in \tint^G$ unless the bond percolation on $G$ given by the set of edges that support at least one bar in $\ladder$ has an infinite component. As noted in \cite{angel}, this implies that 
$\big[ 0 , - \log \big( 1 - p_c(G) \big) \big) \cap \tint^G = \emptyset$, 
where $p_c(G)$ denotes the critical value for bond percolation on $G$. Noting that $p_c(\tree) = \degg^{-1}$, we find that 
\begin{equation}\label{eqlowerbound}
\big[ 0 , \degg^{-1} + \tfrac{1}{2}\degg^{-2} \big) \cap \tint^\tree = \emptyset
\end{equation}
if $\degg \geq 8$.
The papers \cite{angel} and \cite{hammondone} provide two different approaches to proving the existence of infinite cycles in the random stirring process. Appendix $B$ in the arXiv version of~\cite{hammondone} draws on these approaches to provide the following quantitative summary. See Figure~\ref{diffdomains} for an overview of which ranges of $t$ are handled by the two techniques.
\begin{theorem}\label{thmsum}
If $\degg \geq 764$ then $\big[ \degg^{-1} + 2 \degg^{-2}, \infty \big) \subseteq \tint^{\tree}$.
For each $\epsilon > 0$, there exists $d'(\epsilon)$ such that if $\degg \geq d'$ then
 $\big[ \degg^{-1} + ( \tfrac{7}{6} + \epsilon  ) \degg^{-2}, \infty \big) \subseteq \tint^{\tree}$.
\end{theorem}

\begin{figure}
\centering\epsfig{file=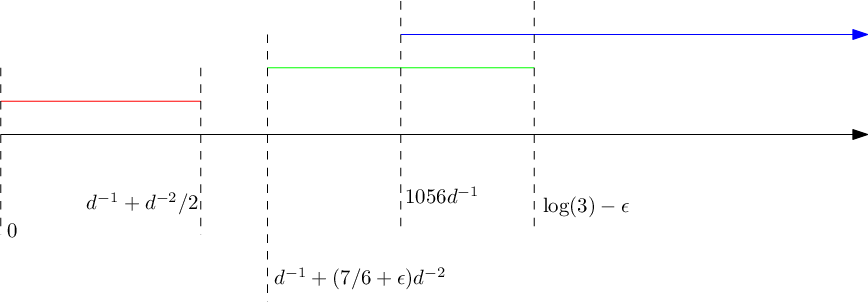, width=14cm}
\caption{The red zone is disjoint from $\tint^{\tree}$ by (\ref{eqlowerbound}). Angel's argument \cite{angel} for infinite cycles works well at small $t$ values, and 
proves that the green zone is contained in $\tint^{\tree}$ if $\degg$ is high. The argument in \cite{hammondone} is valid for all high enough $t$, and shows that the blue zone is  contained in $\tint^{\tree}$ if $\degg$ is high.
}\label{diffdomains}
\end{figure}
\noindent{\bf Proof of Theorem \ref{thm}.}
In light of (\ref{eqlowerbound}) and Theorem~\ref{thmsum}, Theorem~\ref{thm} is reduced to Proposition~\ref{propmono}. \qed
\begin{proposition}\label{propmono}
Suppose that $\degg \geq \numericboundstrong$. Let $\degg^{-1} < s < s' \leq \degg^{-1} + 2\degg^{-2}$. If $s \in \tint^{\tree}$ then $s' \in \tint^{\tree}$. 
\end{proposition}
In fact, our monotonicity result is valid on a slightly longer interval. We now record this other form, since it indicates more accurately the limit of the method; but we will not prove it, since the proof contains no new ideas.
\begin{proposition}\label{propmonorecord}
For any $c \in (0,1)$, there exists $d_0 \in \N$ such that, for $\degg \geq d_0$, if $\degg^{-1} < s < s' \leq (1 - c) \degg^{-1} \log \degg$,
then    $s \in \tint^{\tree}$ implies  $s' \in \tint^{\tree}$. 
\end{proposition}

Cyclic-time random meander $X^\ladder$ is the tool that we will use to prove Propositions~\ref{propmono}. %and~\ref{propmonorecord}. 
Let $\pinf = p_\infty(t)$\hfff{pinf} denote the $\PR_t$-probability that $(\phi,0) \not\in X^\ladder(0,\infty)$. 
Note that in the random stirring model at parameter $t$, the cycle of $\phi$ is infinite 
with probability $p_\infty(t)$.
 
We will prove Proposition~\ref{propmono} 
%(and~\ref{propmonorecord}) 
by establishing that $p_\infty:\big[0, \degg^{-1} + 2\degg^{-2}\big] \to [0,1]$ is non-decreasing (for high enough $\degg$). To do so, we will work with local approximations $\{ p_n: n \in \N \}$ for $p_\infty$.
To define these, we need some notation for describing the graph $\tree$. We pause to collect together such general notation.
\begin{definition} 
 We write $d(\cdot,\cdot): V(\tree) \times V(\tree) \to \N$ for graphical distance on $\tree$.
For $i \in \N$, set $\scv_i = \big\{ v \in V(\tree): d(\phi,v) = i \big\}$\hfff{scv}
and 
$\mathcal{E}_i = \{ e \in E(\tree): d(\phi,e^+) = i \}$.\hfff{edgeset}
For $n \in \N$, let $\treen$\hfff{treen} denote the subtree of $\tree$ induced by vertices at distance at most $n$ from $\phi$.  For $v,w \in V(\tree)$, let $P_{v,w}$\hfff{pathv} denote the unique simple path in $\tree$ connecting $v$ and $w$, and write $E(P_{v,w})$ for its set of edges. 
\end{definition}
Throughout, we write $\vert \cdot \vert$ to denote the cardinality of a set.

For $n \in \N$, let $\reach_n^\ladder \in [0,\infty]$\hfff{hittimen} denote the hitting time $\inf\big\{ s > 0: X^\ladder(s) \in \scv_n \times [0,1) \big\}$. Let $p_n = p_n(t)$\hfff{pn} denote $\PR_t(\reach_n^\ladder < \infty )$. Evidently, $p_n$ decreases pointwise to $\pinf$. Note that $p_0 = 1$.

\begin{subsection}{Pivotality and the added bar}

The tool for deriving Proposition~\ref{propmono}
%and~\ref{propmonorecord} 
is now stated. The derivative of $p_n$ is expressed in terms of the mean effect on $\PR_t(\reach_n^\ladder < \infty)$ caused by adding to $\ladder$ a single ``uniformly'' placed bar. The formula, which is an analogue for the Poisson process of Russo's formula from percolation theory \cite[Theorem 2.25]{grimmett}, is proved in~\cite{zuyev}; however, while Russo's formula is often applied for monotone events in the context of percolation, the event $\big\{ \reach_n^\ladder < \infty \big\}$ is by no means monotone. 
We give a short argument to establish the formula for the sake of completeness. 
\begin{definition} 
Let $n \in \N$. 
Augment the probability space $(\Omega,\PR_t)$
to include a random bar $\addedn$  \hfff{added} 
%= \added_n$ \hfff{added} 
whose law is normalized Lebesgue measure on $E(\treen) \times [0,1)$ and which is independent of $\ladder$. We call $\addedn$ the {\em added bar}. We abuse notation by writing
$\ladder \cup \addedn$ for the bar collection $\ladder \cup \{ \addedn \}$; thus, $X^{\ladder \cup \addedn}$ \hfff{meanderadded} denotes cyclic-time random meander among $\ladder \cup \{ \addedn \}$.
 %We write $\reach_n^\ladder$ and  $\reach_n^{\ladder \cup \addedn}$ for $\reach_n$ for the meanders $X^\ladder$ and  $X^{\ladder \cup \addedn}$.

The on-pivotal event 
%$\onpiv = 
$\onpiv_n$ is defined to be $\big\{ \reach_n^\ladder = \infty,  \reach_n^{\ladder \cup \addedn} < \infty \big\}$,\hfff{onpiv} and the off-pivotal event 
%$\offpiv = 
$\offpiv_n$ to be $\big\{ \reach_n^\ladder < \infty,  \reach_n^{\ladder \cup \addedn} =\infty \big\}$.\hfff{offpiv} 
\end{definition}
\begin{lemma}\label{lemrusso}
For each $n \in \N$,  $p_n:(0,\infty) \to [0,1]$ is differentiable; for $t > 0$, 
$$
 \frac{{\rm d} p_n(t)}{{\rm d} t} = \cardin{E(\treen)} \Big( \PR_t \big( \onpiv_n \big) -   \PR_t \big( \offpiv_n \big) \Big) \, .
$$
\end{lemma}
\noindent{\bf Proof.} Let $\big\{ \ladder_s: s \geq 0 \big\}$ be a coupled system of random bar collections, where $\ladder_s$ has the Poisson-$s$ law on $E(\treen) \times [0,1)$, and where $\ladder_s \subseteq \ladder_{s'}$ whenever $0 \leq s \leq s' < \infty$. Let $N_{s,s'} \in \N$ denote the cardinality of $\ladder_{s'} \setminus \ladder_s$. Note that 
\begin{eqnarray}
 & &  \big\{ \reach_n^{\ladder_t} < \infty \big\}  \cup   \big\{ \reach_n^{\ladder_t} = \infty,  \reach_n^{\ladder_{t+\epsilon}} < \infty,  N_{t,t+ \epsilon}  = 1  \big\}  \cup \big\{  N_{t,t+ \epsilon}  \geq 2 \big\}  \nonumber \\
 & = & \big\{ \reach_n^{\ladder_{t + \epsilon}} < \infty \big\} \cup  \big\{   \reach_n^{\ladder_t} < \infty, \reach_n^{\ladder_{t+\epsilon}} = \infty,  N_{t,t+ \epsilon}  = 1 \big\} \cup \big\{  N_{t,t+ \epsilon}  \geq 2 \big\} \, . \nonumber
\end{eqnarray}
The first two sets in the union of the left-hand side are disjoint. Note that, conditionally on 
$N_{t,t+ \epsilon} = 1$, the unique element in $\ladder_{t + \epsilon} \setminus \ladder_t$ has the distribution of $\addedn$. Thus, taking expectations, we find that
\begin{eqnarray}
 & &p_n(t) +   \epsilon \cardin{E(\treen)}  \exp \big\{ - \epsilon \cardin{E(\treen)} \big\} \PR_t\big(\onpiv_n \big) \nonumber \\
 & \leq & p_n(t+\epsilon) +  \epsilon \cardin{E(\treen)}  \exp \big\{ - \epsilon \cardin{E(\treen)} \big\} \PR_t\big(\offpiv_n \big) +  \big( \epsilon \cardin{E(\treen)} \big)^2 \, . \nonumber 
\end{eqnarray}
which implies that
\begin{equation}\label{eqpnineqone}
p_n(t) +\epsilon \cardin{E(\treen)}   \Big(  \PR_t\big(\onpiv_n \big) -  \PR_t\big(\offpiv_n \big) \Big) \leq p_n(t+\epsilon) + 2 \epsilon^2 \cardin{E(\treen)}^2 \, .  
\end{equation}
Similarly, 
the first two sets in the union of the right-hand side being disjoint, we find that
\begin{equation}\label{eqpnineqtwo}
p_n(t) +\epsilon \cardin{E(\treen)}  \Big(  \PR_t\big(\onpiv_n \big) -  \PR_t\big(\offpiv_n \big) \Big)  + 2 \epsilon^2   \cardin{E(\treen)}^2  \geq p_n(t+\epsilon)\, .  
\end{equation}
From (\ref{eqpnineqone}) and  (\ref{eqpnineqtwo}) follows the statement of the lemma. \qed

\medskip

\begin{proposition}\label{lemloc}
Let $\degg \geq \numericboundstrong$ and suppose that
$\degg^{-1} <t \leq \degg^{-1} + 2\degg^{-2}$. Then, for each $n \geq 1$, $p_n$ is differentiable at $t$, with
$\tfrac{{\rm d} p_n}{{\rm d} t}(t) \geq  \tfrac{\degg}{2} e^{-\degg t} p_n$. 
\end{proposition}

\medskip

\noindent{\bf Proof of  Proposition~\ref{propmono}.} 
By Proposition~\ref{lemloc}, $p_n$ is non-decreasing on $(\degg^{-1},\degg^{-1} + 2\degg^{-2}]$ for each $n \in \N$. 
However, the functions $p_n$ decrease pointwise to $\pinf$.  \qed

\medskip

\subsubsection{The structure of the rest of the paper} We see that the proof of Theorem~\ref{thm} has been reduced to demonstrating Proposition~\ref{lemloc}. 
%(The proof of the incidental Proposition~\ref{propmonorecord} has also been promised, and this will merely involve modifying a few details in the argument for Proposition~\ref{lemloc}.)

In the next section, Section~\ref{secprel}, we present some necessary conditions for on- and off-pivotality, doing so in terms of  the ``crossing'' and ``bottleneck'' events, and we also provide some basic tools. 
We are then in a position in Section~\ref{secmainest} to apply Lemma~\ref{lemrusso} and so reduce Proposition~\ref{lemloc} to two estimates which are there stated: Proposition~\ref{propineqone}, concerning crossing without bottleneck, and Proposition~\ref{propineqtwo}, which treats crossing with bottleneck. 
What remains to complete the proof of Theorem~\ref{thm} is to prove these two estimates. Section~\ref{secfour} provides the proof of Proposition~\ref{propineqone} and Section~\ref{secfive}, that of Proposition~\ref{propineqtwo}.      

%The few details where some change is needed to obtain Proposition~\ref{propmonorecord} occur in the proof of Proposition~\ref{propineqone}, and these changes are explained at the end of Section~\ref{secfour}.

\medskip

Before turning to the preliminaries in Section~\ref{secprel}, we make a comment about the overall approach of the proof of Proposition~\ref{lemloc} via Lemma~\ref{lemrusso}.
The two scenarios depicted in Figure \ref{figcoagfrag} show how the monotonicity
 $\PR_t(\onpiv_n) \geq \PR_t(\offpiv_n)$ is not readily apparent: the appearance of $\addedn$ may lengthen the trajectory of the meander from $(\phi,0)$, so that $\onpiv_n$ occurs, or it may shorten this trajectory and force $\offpiv_n$.  
\begin{figure}
\centering\epsfig{file=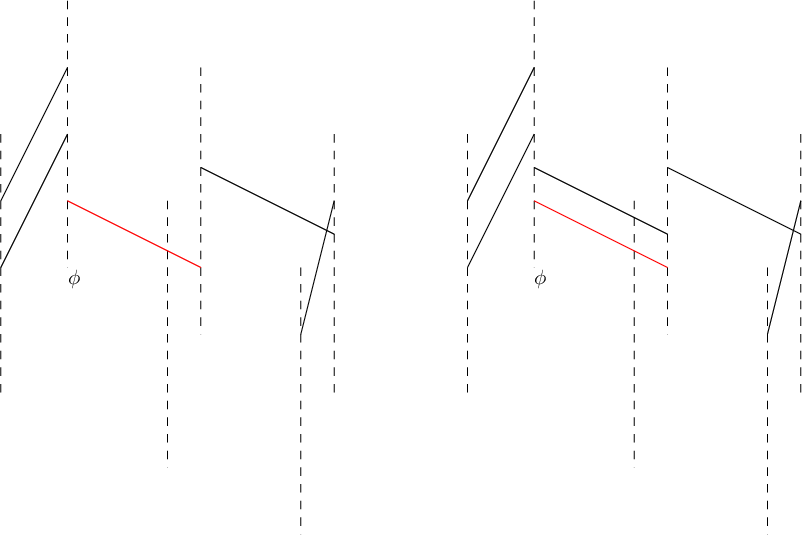, width=10cm}
\caption{In each image, the bars in $\ladder$ are black and the added bar $\addedn$ is red. In the left-hand case, the meander $X^\ladder$ from $(\phi,0)$ completes a circuit in time one; when $\addedn$ appears, the meander $X^{\ladder \cup \addedn}$ makes a longer journey, perhaps never returning to $(\phi,0)$. In the right-hand case, the appearance of $\addedn$ has the opposite effect, curtailing the trajectory of the meander from~$(\phi,0)$.}\label{figcoagfrag} 
\end{figure}
To suggest our approach in a few words, we will argue that the coagulating mechanism causing $\onpiv_n$  is stronger than the fragmenting one causing $\offpiv_n$ when $t = \Theta(\degg^{-1})$ because, for such $t$, the bar collection $\ladder$ is dilute: the added bar $\addedn$ will probably arrive over an edge where no bar of $\ladder$ is present, and then (as we will prove shortly in Lemma~\ref{lemdoubonethree})  
the meander $X^{\ladder \cup \addedn}$ follows either the same route as does $X^\ladder$ or a longer one. 

\medskip 

\noindent{\bf Acknowledgments}. I would like to thank two referees, one for  suggesting a new guise for the argument which has led to a significant simplification in Section~\ref{secfour}, and the other for providing a thorough critique of an  earlier version.

\end{subsection}
\end{section}
\begin{section}{Conditions for pivotality and basic tools}\label{secprel}

\begin{subsection}{Three necessary conditions for pivotality of the added bar}

We will introduce several events in order to discuss this pivotality: the crossing, bottleneck and non-escape events, $\cross_n,\seev_n$ and $\nonesc_n$. To keep formulas uncluttered, we will omit the $n$-dependence in denoting these events, as we will also for the pivotality events $\onpiv_n$ and $\offpiv_n$. 

\begin{subsubsection}{The meander must encounter the added bar}

Let $\cross = \cross_n$ \hfff{cross} denote the {\it crossing} event that $X^\ladder$ meets a joint of $\addedn$ before time $\reach_n^\ladder$. If $\cross$ does not occur, then the trajectories of $X^\ladder$ and $X^{\ladder \cup \addedn}$ are equal at least on the interval $[0,\reach_n]$ (where the value of $\reach_n$ is shared by the two processes); this proves the following fact.
\begin{lemma}\label{lemcross}
We have that
$$
\onpiv \cup \offpiv \subseteq \cross \, .
$$
\end{lemma}
\end{subsubsection}
\begin{subsubsection}{No escape above the bottleneck bar}
 
 In the case that $\cross$ occurs, we now provide a further necessary condition for the occurrence of 
$\onpiv \cup \offpiv$. 
% For $v,w \in V(\tree)$, let $P_{v,w}$\hfff{pathv} denote the unique simple path in $\tree$ connecting $v$ and $w$, and $E(P_{v,w})$ for its set of edges. 
If $\cross$ occurs, note that each element of $E\big(P_{\phi,E(\addedn)^+}\big)$ supports at least one bar in $\ladder$. The bottleneck event $\seev$\hfff{seev} occurs if one of these elements supports exactly one bar in $\ladder$.  
If $\seev$ occurs,
define the {\it bottleneck edge} $\singleedge$ \hfff{singleedge} to be the edge on $P_{\phi,E(\addedn)^+}$  supporting exactly one bar in $\ladder$ that is furthest from the root. Let $\singlebar$ \hfff{singlebar} denote the unique bar on $\singleedge$.
If $\cross \cap \seev$ occurs,
then $X^\ladder$ certainly crosses $\singlebar$.  
If also $X^\ladder$ has a periodic trajectory, then $X^\ladder$ must later cross back along 
$\singlebar$ to arrive at $\singlebar^+$.
The  {\it non-escape} event $\nonesc$ \hfff{nonesc} occurs if the 
meander  $X_{\singlebar^+}^\ladder$ visits $(\phi,0)$ before $\scv_n \times [0,1)$. 
We claim that
\begin{equation}\label{eqnonesc}
\cross \cap \seev \cap \{ \reach_n^\ladder = \infty \} \subseteq \nonesc \, . 
\end{equation}
Indeed, as we have seen, occurrence of
 the left-hand event implies that  $X^\ladder$ at some time recrosses $\singlebar$ to arrive at  $\singlebar^+$; after this time, $X^\ladder$
follows the route of $X_{\singlebar^+}^\ladder$, so that $\reach_n^\ladder = \infty$ forces $\nonesc$, and we have~(\ref{eqnonesc}).
The inclusion~(\ref{eqnonesc}) holds equally if $\reach_n^{\ladder \cup \addedn}$ replaces $\reach_n^\ladder$; the same argument works after we note that $\addedn$ is supported on an edge in the descendent tree of $\singleedge^-$,  %E\big(\tree_{[\singleedge^-]}\big) \times [0,1)$, 
and thus $X^{\ladder \cup \addedn}_{\singlebar^+}$ and $X^{\ladder}_{\singlebar^+}$ coincide at least until return to 
$\singlebar^+$, by which time the two processes have visited $(\phi,0)$ because $X^\ladder$ (from $(\phi,0)$) visits $\singlebar^+$.
These inferences form the basis for the following claim.
\begin{lemma}\label{lemaddededge}
We have that
$$
\big( \onpiv \cup \offpiv \big) \cap \seev \subseteq \cross  \cap \nonesc \, .
$$
\end{lemma}
\noindent{\bf Proof.}
It suffices in light of Lemma \ref{lemcross} to argue that $\cross \cap \seev \cap \nonesc^c \cap \big( \onpiv \cup \offpiv \big) = \emptyset$. However, we have argued that 
$\cross \cap \seev \cap \nonesc^c$ 
forces both $\reach_n^\ladder < \infty$ and $\reach_n^{\ladder \cup \addedn} < \infty$; this suffices,  
because then neither $\onpiv$ nor $\offpiv$ may occur. \qed
\end{subsubsection}
\begin{subsubsection}{For off-pivotality, the edge of the added bar must support a bar in $\ladder$}
\begin{lemma}\label{lemdoubonethree}
If $E(\addedn)$ supports no bar in $\ladder$ then $\offpiv$ does not occur.
\end{lemma}
For use in the proof, write $\hit_{\addedn} = \inf \big\{ s> 0: X(s) \in \{ \addedn^+,\addedn^- \} \big\}$,\hfff{hitadded} and note that $X^{\ladder}$ and $X^{\ladder \cup \addedn}$ coincide until $\hit_{\addedn}$, whose value the two processes share.

\medskip

\noindent{\bf Proof of Lemma~\ref{lemdoubonethree}.}
By Lemma \ref{lemcross}, we may assume that $\cross$ occurs.
If~$E(\addedn)$ supports no bar in $\ladder$, then the trajectory $X^{\ladder \cup \addedn}$ is formed from that of $X^\ladder$ as follows. Note that the two processes reach the parent joint of $\addedn$ at time $\hit_{\addedn}$: it is impossible that $X^\ladder$ be at the pole of $\addedn^-$ at any time, because this would entail $X^\ladder$ crossing $E(\addedn)$, an edge which supports no bar in $\ladder$. Thus, at time $\hit_{\addedn}$ it is from the parent to the child joint that $X^{\ladder \cup \addedn}$ crosses $\addedn$. This meander then spends a duration 
among the poles associated to
the descendent tree of $E(\addedn)^-$.
%in $E\big(\tree_{[E(\addedn)^-]}\big) \times [0,1)$. 
It may visit $\scv_n \times [0,1)$ during this sojourn, so that $\reach_n^{\ladder \cup \addedn} < \infty$ occurs, thus excluding~$\offpiv$. If this does not happen, $X^{\ladder \cup \addedn}$ recrosses $\addedn$ to reach $\addedn^+$ again and then continues to follow the trajectory of $X^\ladder$ from this point. The two processes have no further opportunity to diverge (except by a return to $\addedn^+$ in a later circuit), and this makes $\offpiv$ impossible. \qed
\end{subsubsection}
\end{subsection}
\begin{subsection}{Some basic tools}
Here we record two simple observations regarding cyclic-time random meander.
\begin{lemma}\label{lemsonesym}
The distribution of the return time to $(\phi,h)$ of $X^\ladder_{(\phi,h)}:[0,\infty) \to V(\tree) \times [0,1)$ under $\PR_t$ is independent of $h \in [0,1)$. 
\end{lemma}
\noindent{\bf Proof.}
The bar collection $\ladder$ has the Poisson-$t$ distribution on $E(\tree) \times [0,1)$ and thus is invariant under the map which increases the height of all bars by $h$ and reduces modulo $1$. \qed
\begin{lemma}\label{lemubl}
Let $s > 0$.
 Consider the law $\PR_t$ given $X^\ladder:[0,s] \to V(\tree) \times [0,1)$. Let $\forced_s \subseteq E(\tree) \times [0,1)$
%\hfff{forcedname} 
denote the set of bars in $\ladder$ that $X^\ladder$ has crossed during $[0,s]$. Let the set of time-$s$ {\em untouched} bar locations $\ubl_s \subseteq E(\tree) \times [0,1)$ 
%\hfff{untouched} 
denote the set of bars $b \in E(\tree) \times [0,1)$ neither of whose joints belongs to $X[0,s]$. Then the conditional distribution of $\ladder$ is given by $\forced_s \cup \ladder_{(s,\infty)}$, where $\ladder_{(s,\infty)}$ is a random bar collection with Poisson law of density $t 1\!\!1_{\ubl_s}$.
\end{lemma}
\noindent{\bf Proof.} That $\forced_s \subseteq \ladder$ is known given $X^\ladder$ on $[0,s]$; similarly, if $X^\ladder[0,s]$ visits the joint of some bar in $\ladder$, 
that bar belongs to $\forced_s$. The time-$0$ distribution of the remaining bars, those in $\ubl_s$, is undisturbed by the data~$X^\ladder[0,s]$. \qed 

\end{subsection}
\end{section}
\begin{section}{Deriving monotonicity: the proof of  Proposition~\ref{lemloc}}\label{secmainest}

By Lemma \ref{lemcross}, in comparing the probabilities of $\onpiv$ and $\offpiv$, we may restrict attention to choices of $\ladder$ and $\addedn$ such that $\cross$ occurs. We will further divide into cases according to whether $\seev$ occurs, and, if it does, according to the location of the bottleneck edge $\singleedge$.

We now record how we split into such cases, state the estimates that we will prove for each case, and then provide the proof of Proposition~\ref{lemloc} using these estimates.

Lemmas \ref{lemcross} and \ref{lemaddededge}
imply that the right-hand side of the following equality is a partition into disjoint sets:  for $* \in \{ +, - \}$,
$$
\mathsf{P}^*   = 
\mathsf{P}^*   \cap \cross \cap \seev^c 
 \, \bigcup  \, 
\mathsf{P}^*   \cap \cross  \cap \seev \cap \nonesc   \, .
$$
We have  then that
\begin{equation}\label{eqthreea}
 \PR_t \big( \onpiv \big) -  \PR_t \big( \offpiv \big) = \nta_1 + \nta_2  \, ,
\end{equation}
where 
\begin{equation}\label{eqaone}
\nta_1 = 
 \PR_t \big( \onpiv  \cap \cross  \cap \nse \big) -  \PR_t \big( \offpiv  \cap \cross   \cap \nse \big) 
\end{equation}
and
$$
\nta_2 = 
 \PR_t \big( \onpiv   \cap \cross  \cap \seev \cap \nonesc  \big) -  \PR_t \big( \offpiv  \cap \cross  \cap \seev  \cap \nonesc  \big) \, .
$$
Note that, under $\seev$, $\singleedge \in E(\tree_{n-1})$, and thus
\begin{equation}\label{eqatwo}
\nta_2 = \sum_{e \in E(\tree_{n-1})} \Big(
 \PR_t \big( \onpiv  \cap \cross \cap \nonesc \cap \big\{ \singleedge = e \big\} \big) \, - \,  \PR_t \big( \offpiv  \cap \cross \cap \nonesc \cap \big\{ \singleedge = e \big\} \big) \Big) \, .
\end{equation}
\begin{subsection}{The main estimates}

We set $\tau = t \degg$,\hfff{tau} and will often use $\tau$ in place of $t$ in the ensuing estimates. Note that the range of $t$-values that concerns us, 
$[\degg^{-1},\degg^{-1} + 2\degg^{-2}]$, corresponds to $\tau \in [1,1 + 2\degg^{-1}]$, so that $\tau$ is a unit-order quantity. Our use of the notation $\tau$ will emphasise how our argument is perturbative in high $\degg$:
for example, we will see that, up to $\tau$-dependent unit-order factors, $\PR_t \big( \onpiv  \cap \cross  \cap \nse \big)$ is at least $\degg  \tfrac{p_{n-1}}{\cardin{E(\treen)}}$, and  $\PR_t \big( \offpiv  \cap \cross   \cap \nse \big)$ is at most   $\tfrac{p_{n-1}}{\cardin{E(\treen)}}$, so that $\nta_1$ in~(\ref{eqaone}) is positive if $\degg$ is high enough.

We now state the two main estimates that will lead to Proposition~\ref{lemloc}. 

\begin{proposition}\label{propineqone} 
Let $n \geq 1$. If  $\degg \geq \numericboundstrong$  and $\tau \in \big[1 , 1 +  2\degg^{-1} \big]$
then 
$$
  \PR_t\big(\onpiv \cap \cross \cap \seev^c \big) \geq   \PR_t\big(\offpiv \cap \cross \cap \seev^c \big) \, + \, 
   \tfrac{\degg}{2} e^{-\tau} \tfrac{p_{n-1}}{\cardin{E(\treen)}}  \, .
$$
\end{proposition}
\begin{proposition}\label{propineqtwo}
Suppose that $n \geq 2$, $\degg \geq \numericboundstrong$ and $\tau \in \big[1 , 1 +  2\degg^{-1} \big]$.
Then, for all $e \in E(\tree_{n-1})$,
$$
 \PR_t \Big( \onpiv \cap \cross \cap \nonesc \cap \big\{ \singleedge = e \big\}  \Big) \geq \PR_t \Big( \offpiv  \cap \cross \cap \nonesc \cap  \big\{ \singleedge = e \big\}  \Big) \, .
 $$
\end{proposition}

\end{subsection}
\begin{subsection}{Combining the estimates}

We now apply the estimates in Propositions~\ref{propineqone} and~\ref{propineqtwo}  
to Lemma~\ref{lemrusso} in order to prove  Proposition~\ref{lemloc}. 

\medskip

\noindent{\bf Proof of Proposition~\ref{lemloc}.}
By Propositions~\ref{propineqone} and~\ref{propineqtwo}, and~(\ref{eqaone}) and~(\ref{eqatwo}),
\begin{equation}\label{eqaonetwo}
 A_1 + A_2 \geq \tfrac{\degg}{2} e^{-\tau} \tfrac{p_{n-1}}{\cardin{E(\treen)}} \,.
\end{equation}

The proposition then follows from Lemma~\ref{lemrusso},~(\ref{eqthreea}) and $p_{n-1} \geq p_n$.  \qed

 \medskip
 
It remains to prove Propositions \ref{propineqone} and~\ref{propineqtwo}. Sections~\ref{secfour} and~\ref{secfive} treat each of these estimates in turn.
\end{subsection}
\end{section}
\begin{section}{Crossing without bottleneck: the proof of Proposition~\ref{propineqone}}\label{secfour}

Some definitions are needed for this proof.
Say that an element $e \in E(\tree)$ is multi-open if $e$ supports at least two bars in $\ladder$. 
For $v \in V(\tree)$, the  multi-open component of $v$ is the set of $w \in V(\tree)$ such that every edge in $E(P_{v,w})$ is multi-open.  
Define the multi-cluster $\double_v$\hfff{multicl} to be the edge-set of the subgraph of $\tree$ induced by the multi-open component of $v$.

Given $G \subseteq E(\tree)$, let $\pext G$ \hfff{pext} denote the set of edges in $E(\tree) \setminus G$
that are incident to the endpoint of some element of $G$.

\begin{lemma}\label{lemdoublesize}
Let $\degg \geq 1$. For each $\ell \in \N^+$,
$$
\PR_t \big( \cardin{\doublephi} = \ell  \big) \leq   e^{-1} \big( e \tau^2 \degg^{-1} \big)^\ell \, .
$$
\end{lemma}
\noindent{\bf Proof.}
Consider the following procedure for determining $\doublephi$, which is similar to the coding of trees by Lukasiewicz paths presented in \cite[Section 1.1]{legall}. 
At time zero, the $\degg$  
elements of $\mathcal{E}_0$ are called candidates; at each time step from time one, one candidate is examined. 
On being examined, a candidate changes status, either being found to belong to $\doublephi$, or not. In the first case, the edges incident to the candidate's child vertex join the candidate list; in the second, no new candidates join. The process stops when there are no candidates left. The set of edges that are candidates at some time is 
$\doublephi \cup \pext \doublephi$. 

For $r \in [0,1]$, let $Z = Z_r:\N \to \Z$, $Z(0) = \degg$, denote the Markov chain on $\Z$ with two transitions,
namely a $\degg-1$ displacement with probability $r$ and a $-1$ displacement with probability $1-r$.
Let $\sigma_r = \inf \big\{ \ell \geq 0 : Z(\ell) = 0 \big\}$. 

Henceforth we write $Z$ and $\sigma$ for $Z_r$ and $\sigma_r$ with $r = 1 - (1 + t)e^{ - t }$.
Any given candidate belongs to $\doublephi$ with probability $1 - (1 + t)e^{ - t }$.
Hence, for each $s \in \big\{ 0,\ldots, \sigma \big\}$, $Z(s)$ is the number of candidates after the examination at time $s$. We see that 
$\cardin{\doublephi}  + \cardin{\pext \doublephi}$ under $\PR_t$
 has the distribution of 
$\sigma$.
Noting that $\cardin{\doublephi} +  \cardin{\pext \doublephi} = \degg \big( \cardin{\doublephi} + 1 \big)$, we find that
$$
  \PR_t \big( \cardin{\doublephi}  = \ell   \big) =   \PR_t \big(  \cardin{\doublephi}  +  \cardin{\pext \doublephi}  = \degg (\ell + 1)  \big) 
 = \PR \big( \sigma = \degg (\ell +1)  \big) \, .  
$$
Note now that since the occurrence of $\sigma = \degg (\ell +1)$ implies that $Z\big( \degg(\ell + 1) \big) = 0$,
it entails that  exactly $\ell$ among the first $\degg (\ell+1)$ transitions made by $Z$ are up moves; moreover, 
since this event also requires that $Z$ has no earlier visit to zero, each of the last $\degg$ such transitions are down moves, so that, in fact, exactly $\ell$ among the first $\degg \ell$ transitions of $Z$ are up moves. Also using $1 - (1 + t)e^{ - t } \leq t^2$, we find that 
$$
\PR_t \big( \cardin{\doublephi}  = \ell   \big) \leq  {\degg \ell \choose \ell} t^{2\ell} \, ,
$$
whose right-hand side is at most 
 $\tfrac{(\degg \ell)^\ell}{\ell !} t^{2 \ell}$ which, since $\ell \geq 1$, is bounded above by $\degg^\ell e^{\ell - 1} t^{2\ell} = e^{-1} \degg^{-\ell} e^\ell \tau^{2\ell}$.  \qed

Let $\single_\phi$ \hfff{singlephi} denote the set of elements in $\pext \double_\phi$ that support a bar in $\ladder$. Note that each edge in $\single_\phi$ supports exactly one bar. 

\begin{lemma}\label{lempoisson}
Let $t > 0$ and $k \in \N$. Under the law $\PR_t$ given that $\cardin{\doublephi} = k$, the conditional distribution of  $\cardin{\singlephi}$
is stochastically dominated by a Poisson random variable of parameter $(k+1)\tau$.  
\end{lemma}
\noindent{\bf Proof.} The random variable $\cardin{\singlephi}$ is a sum over $\pext \doublephi$ of independent Bernoulli random variables. Each of these random variables has the law of a Poisson random variable of parameter $\tau \degg^{-1}$ conditioned to assume value either zero or one, a law which is stochastically dominated by that of a Poisson random variable of parameter $\tau \degg^{-1}$.
Thus, conditionally on $\vert \pext \doublephi \vert$,  $\cardin{\singlephi}$ is stochastically dominated by a Poisson random variable of parameter $\vert \pext \doublephi \vert \tau \degg^{-1}$. 
The result follows from $\vert \pext \doublephi \vert = \degg + (\degg - 1) \cardin{\doublephi}$. \qed

\begin{lemma}\label{lempcblb}
Suppose that $\degg \geq 2$, $n \geq 1$ and $t > 0$. Then
$$
  \PR_t\big(\onpiv \cap \cross \cap \seev^c \big)  \geq \degg e^{-\tau}  \tfrac{p_{n-1}}{\cardin{E(\treen)}} \, .
$$
\end{lemma}
\noindent{\bf Proof.} The event that $\mathcal{E}_0 \times [0,1)$ is disjoint from $\ladder$ but contains $\addedn$ has $\PR_t$-probability $\degg  \big( \cardin{E(\treen)} \big)^{-1} e^{-\tau}$; and, when it occurs,  so does
$\cross \cap \seev^c$. Moreover, when it occurs,
 clearly $\hit_n^\ladder = \infty$, while $\hit_n^{\ladder \cup \addedn} < \infty$ precisely when the meander $X^{\ladder \cup \addedn}_{\addedn^-}$ visits $\scv_n \times [0,1)$ before its return to $\addedn^-$, an independent event of probability $p_{n-1}$. \qed

\begin{lemma}\label{lempncomp}
Let $\degg \geq 2$ and $t \geq 0$. For $n,m \in \N$ such that $n \geq m$, $p_n \geq \big( ( 1 - e^{-\degg t} ) e^{-t} \big)^{n-m} p_m$. 
\end{lemma}
\noindent{\bf Proof.}
We have that $p_1 = 1 - e^{-\degg t}$, because $X^\ladder$ leaves the pole of $\phi$
if and only if an edge incident to $\phi$ supports a bar in $\ladder$. In view of this, it is enough to
establish that, for any $n > m \geq 1$,
\begin{equation}\label{e.pnpm}
p_n \geq e^{-t} p_m p_{n-m} \, .
\end{equation}

To see this, note that, given $H_m^\ladder < \infty$,  $X^\ladder$ crosses a bar from its parent to its child joint at time $H_m^\ladder$. Call this bar $b$ and the edge which supports it $e$, and note that $d(\phi,e^-) = m$. There is, by Lemmas~\ref{lemsonesym} and~\ref{lemubl}, conditional probability at least $e^{-t} p_{n-m}$ that $e$ supports no bar but $b$, and that the vertex component of $X^\ladder_{b^-}$ reaches distance $n-m$ from $e^-$ before returning to $e^-$. These circumstances force the occurrence of $H_n^\ladder < \infty$. Thus, we obtain~(\ref{e.pnpm}) and so complete the proof. \qed

\begin{lemma}\label{lempcbub}
Suppose that $n \geq 1$, $\degg \geq \tfrac{32 e^3}{e-1}$ and $1 \leq \tau \leq 2$. Then
$$
  \PR_t \big( \offpiv \cap \cross \cap \seev^c \big)  \leq
   \tfrac{p_{n-1}}{\cardin{E(\treen)}} \Big(   \tau(\tau + 1)
   \, + \,   \tfrac{768 e^2}{e-1}  \degg^{-1} \Big) \, .
$$
\end{lemma}
\noindent{\bf Proof.} By definition, $\offpiv \subseteq \{ \hit_n^\ladder < \infty \}$. Note also that $\cross \cap \seev^c$ entails that $E(\addedn) \in \doublephi \cup \pext \doublephi$; if also $\offpiv$ occurs, then, by Lemma~\ref{lemdoubonethree}, 
$E(\addedn) \in \pext \doublephi$ is possible only if $E(\addedn) \in \singlephi$. Hence, 
\begin{equation}\label{eqoffpivinc}
\offpiv \cap \cross \cap \seev^c \subseteq \cup_{k,s \geq 0}  
  \big\{ \hit_n^\ladder < \infty \big\} \cap \big\{ \cardin{\doublephi} = k , \cardin{\singlephi} = s ,
  E(\addedn) \in \doublephi \cup \singlephi \big\} \, .
\end{equation}

We now claim that, for $k,s \in \N$ such that $n \geq k+1$, 
\begin{equation}\label{e.hitnbd}
 \PR_t \Big( \hit_n^\ladder < \infty \, \Big\vert \, \cardin{\doublephi} = k , \cardin{\singlephi} = s  \Big) \leq s p_{n - k - 1} \, .
\end{equation}

To derive this, consider the event that
$\cardin{\doublephi} = k$ and $\cardin{\singlephi} = s$. 
Let
$B_1,\ldots,B_s$ denote the set of bars in $\ladder$ supported on elements of $\singlephi$. 
Note that each child joint $B_i^-$ lies on a pole whose vertex is at distance at most $\cardin{\doublephi} + 1 \leq k + 1$ from $\phi$, so that $B_i^- \in V(\tree_{n}) \times [0,1)$ by assumption. If  $\hit_n^\ladder < \infty$ is to occur, then, for some $1 \leq i \leq s$, it is necessary that the meander $X_{B_i^-}^\ladder$ reaches $\scv_n \times [0,1)$ before its return to $B_i^-$. For given $i \in \N$, under the law $\PR_t$ conditioned on any given instance of the intersection of $\ladder$ and $\big( \doublephi \cup \singlephi \big) \times [0,1)$ (for which $\cardin{\singlephi} \geq i$), the conditional probability that  $X_{B_i^-}^\ladder$ does so equals $p_{n - d \big( \phi, V(B_i^-) \big)} \leq p_{n - k - 1}$. Thus we obtain~(\ref{e.hitnbd}).

By (\ref{e.hitnbd}), the form of the law of $\addedn$ and the independence of $\addedn$ and $\ladder$, we find that, if $n \geq k+1$,
\begin{eqnarray}
 & & \PR_t \Big( \hit_n^\ladder < \infty , \cardin{\doublephi} = k , \cardin{\singlephi} = s , E(\addedn) \in \doublephi \cup \singlephi \Big) \label{eqhitnbd} \\
 & \leq & 
\tfrac{k + s}{\cardin{E(\treen)}} \PR_t \big(  \cardin{\doublephi} = k , \cardin{\singlephi} = s  \big) s p_{n - k -1}  \, . \nonumber
\end{eqnarray}
%Adopting the convention that $p_i = 1$ for $i \leq 0$, this inequality may be extended to hold for all $n \in \N$ when it takes the form 

When $n \leq k$, we have the bound
\begin{eqnarray}
 & & \PR_t \Big( \hit_n^\ladder < \infty , \cardin{\doublephi} = k , \cardin{\singlephi} = s , E(\addedn) \in \doublephi \cup \singlephi \Big) \label{eqhitnbdnew} \\
 & \leq & 
\tfrac{k + s}{\cardin{E(\treen)}} \PR_t \big(  \cardin{\doublephi} = k , \cardin{\singlephi} = s  \big)  \, . \nonumber
\end{eqnarray}

Adopt the convention that $p_i = 1$ for $i < 0$ (as well as for $i = 0$, as we already prescribed). 
For $n \in \N$, define 
$$
 A_{n,0} = \sum_{s=0}^\infty \tfrac{s^2}{\cardin{E(\treen)}} \PR_t \big(  \cardin{\doublephi} = 0 , \cardin{\singlephi} = s  \big) p_{n  - 1}  \, ,
$$
and, for $k \geq 1$, define
$$
 A_{n,k} = \sum_{s=0}^\infty \tfrac{(k + s)(s+1)}{\cardin{E(\treen)}} \PR_t \big(  \cardin{\doublephi} = k , \cardin{\singlephi} = s  \big) p_{n - k -1}  \, ,
$$

Note from~(\ref{eqoffpivinc}), (\ref{eqhitnbd}) and (\ref{eqhitnbdnew}) that
\begin{equation}\label{eqaksum}
 \PR_t \big( \offpiv \cap \cross \cap \seev^c \big) \leq \sum_{k=0}^\infty \akn \, .
\end{equation}
We have that
\begin{eqnarray}
 \aon & \leq &  \tfrac{p_{n-1}}{\cardin{E(\treen)}} \sum_{s=0}^{\degg} s^2  \PR_t \big(   \cardin{\singlephi} = s
 \, \big\vert \, \cardin{\doublephi} = 0
 \big) \nonumber \\
  & \leq &    \tfrac{p_{n-1}}{\cardin{E(\treen)}} \mathbb{E} Z_\tau^2 =   \tfrac{p_{n-1}}{\cardin{E(\treen)}} \tau \big( \tau + 1 \big) \, , \label{eqakzero}
\end{eqnarray}
where here $Z_s$, for $s \geq 0$, has the Poisson distribution of parameter $s$. The second inequality above invokes Lemma~\ref{lempoisson}.

When $k \geq 1$, we find from Lemma~\ref{lemdoublesize} that
$$
\akn \leq  \tfrac{p_{n-k-1}}{\cardin{E(\treen)}}   e^{-1} \big( e \tau^2 \degg^{-1} \big)^k \sum_{s=0}^{k(\degg - 1) + \degg} (k + s)(s+1) \, \PR_t \big(   \cardin{\singlephi} = s
 \, \big\vert \, \cardin{\doublephi} = k 
 \big) \, .
$$
Note that, by Lemma~\ref{lempoisson},
\begin{eqnarray*}
 & & \sum_{s=0}^{k(\degg - 1) + \degg} (k + s)(s+1) \, \PR_t \big(   \cardin{\singlephi} = s
 \, \big\vert \, \cardin{\doublephi} = k 
 \big) \\
 &\leq & (k + 1) \mathbb{E} Z_{(k+1)\tau} +  \mathbb{E} Z_{(k+1)\tau}^2 + k = (k+1)^2\tau + (k+1)^2 \tau^2 + (k + 1)\tau + k \, . 
\end{eqnarray*}
Thus, for $k \geq 1$,
\begin{equation}\label{eqaksummary}
\akn \leq  \tfrac{p_{n-k-1}}{\cardin{E(\treen)}}  e^{-1} \big( e \tau^2 \degg^{-1} \big)^k  \big( (k+1)^2 \tau (\tau + 1) + (k+1) \tau   + k \big) \, .
\end{equation}

Lemma~\ref{lempncomp} bounds above $\tfrac{p_{n-k-1}}{p_{n-1}}$ when $n- k -1 \geq 0$, and, given our convention that $p_i = 1$ for $i < 0$, this bound may trivially be extended to hold for all values of $n-k-1$. Using the bound,  and $\tau \leq 2$, we find that
$$
\akn \leq  \tfrac{p_{n-1}}{\cardin{E(\treen)}}    e^{-1} \big( (1 - e^{-\tau})^{-1} e^{\tau/\degg} e \tau^2 \degg^{-1} \big)^{k} \cdot 6 (k+2)^2 \, .
$$
By $1 \leq \tau \leq 2$, $\degg \geq 2$ and $\tau^2 \degg^{-1} \leq 1/2$, 
$$
\akn \leq  \tfrac{p_{n-1}}{\cardin{E(\treen)}}  6  e^{-1} \big( \tfrac{4 e^3}{(e-1) \degg} \big)^{k} (k+2)^2  \, .
$$

Using $k+2 \leq 2^{k+1}$ and $\degg \geq \tfrac{32 e^3}{e-1}$,
\begin{equation}\label{eqakone}
\sum_{k=1}^\infty \akn \leq   \tfrac{768 e^2}{e-1} \tfrac{p_{n-1}}{\cardin{E(\treen)}}     \degg^{-1}    \, .
\end{equation}
We obtain Lemma~\ref{lempcbub} from~(\ref{eqaksum}),~(\ref{eqakzero}) and~(\ref{eqakone}). \qed

\medskip

We may now combine Lemmas~\ref{lempcblb} and~\ref{lempcbub} to prove Proposition~\ref{propineqone}.

\medskip

\noindent{\bf Proof of Proposition~\ref{propineqone}.}
Recall from the proposition's statement that we are assuming that $\degg \geq \numericboundstrong$ and that $1 \leq \tau \leq 1 + 2 \degg^{-1}$, conditions which imply that each of the following is satisfied: 
$\degg \geq \tfrac{32 e^3}{e-1}$, $1 \leq \tau \leq 2$ and  
$$
 \tfrac{\degg}{2} e^{-\tau} \geq   \tau(1+\tau) + 
   \tfrac{768 e^2}{e-1} \degg^{-1}  \, ,
$$    

Thus, the proposition follows from  Lemma~\ref{lempcblb} and  Lemma~\ref{lempcbub}. \qed

\end{section}
\begin{section}{Crossing with bottleneck: the proof of Proposition~\ref{propineqtwo}}\label{secfive}

In this section, we prove Proposition~\ref{propineqtwo}. The argument when formally recorded has a technical appearance but in fact, as we now explain, it is a straightforward reduction to the crossing without bottleneck  Proposition~\ref{propineqone}. First a definition is convenient.
\begin{definition}\label{def.desc}
For $v \in V(\tree)$, let $\desctree{v}$  \hfff{desctree}  denote the subtree of $\tree$ induced by descendents of $v$ ($\desctree{v}$ may viewed as a rooted tree with root $v$).
\end{definition}

Proposition~\ref{propineqtwo} will be proved by a study of the law $\PR_t \big(\cdot \big\vert \cross \cap \seev \cap \nonesc \big)$. Under this measure, $X^{\ladder \cup \addedn}$ will necessarily cross the bottleneck bar to arrive at its child joint $\singlebar^-$. It then makes an {\em excursion} over the descendent tree $\desctree{\singleedge^-}$, following the trajectory of 
the process $X_{\singlebar^-}$ until returning, at a time which is perhaps infinite, to $\singlebar^-$; moreover, in order that $\cross$ be realized, the  process must, during its excursion, cross $\addedn$ before it may either return to $\singlebar^-$ or reach the boundary $\scv_n \times [0,1)$; and, in order that the identity of the bottleneck bar $\singlebar$ be respected, the process, from the start of the excursion until its crossing of $\added_n$, may cross only over edges supporting at least two bars of $\ladder$. Loosely put, the process during the excursion under the conditioned measure verifies the crossing without bottleneck event, where this event is associated to the tree $\desctree{\singleedge^-}$
and starting point $\singlebar^-$ in place of the whole tree and $(\phi,0)$. More precisely, 
a moment's thought shows that (and, in our formal verification below, Lemma~\ref{lemseev} states that), under our conditioning, the joint law of $\big( \ladder \cup \added_n \big) \cap \big( \desctree{\singleedge^-} \times [0,1) \big)$ and $X^{\ladder \cup \added_n}$ during the excursion has the following form.  The root  is taken to be the vertex $\singleedge^-$. The quantity $n$ is replaced by  $n - d(\phi,\singleedge^-)$, so that it still measures the distance from the root to the boundary. Time is cyclically shifted so that $\singlebar^-$ has height zero (so that, in the new coordinates, the process at the start of excursion is at the root and at height zero). Then, after this relabelling, the joint law above is simply given by conditioning $\PR_t$ on $\cross_{n - d(\phi,\singleedge^-)} \cap \seev_{n - d(\phi,\singleedge^-)}$, that is, by conditioning on the crossing without bottleneck event when the added bar $\added_{n - d(\phi,\singleedge^-)}$ appears uniformly on $E \big( \tree_{n - d(\phi,\singleedge^-)} \big) \times [0,1)$. 

\medskip

How does this understanding of the conditioning on crossing with bottleneck allow us to prove  Proposition~\ref{propineqtwo}? Clearly, to prove the proposition, it is enough to argue that,  under its hypotheses, for each $(e,h) \in E(\tree_{n-1}) \times [0,1)$, 
\begin{equation}\label{eqprtsb}
 \PR_t \big( \onpiv_n \big\vert \cross_n \cap \nonescn \cap \big\{ \singlebarn = (e,h) \big\}  \big) \geq 
\PR_t \big( \offpiv_n  \big\vert \cross_n \cap \nonescn \cap  \big\{ \singlebarn = (e,h) \big\}  \big) \, .
\end{equation}

Under the conditioning in (\ref{eqprtsb}), and in the new coordinates just described, the events $\onpiv_n$ and $\offpiv_n$ translate to $\onpiv_{n-d(\phi,\singleedge^-)}$ and $\offpiv_{n-d(\phi,\singleedge^-)}$. See Lemma~\ref{lemseevplus} below. Since the conditioning translates to conditioning on $\cross_{n - d(\phi,\singleedge^-)} \cap \seev_{n - d(\phi,\singleedge^-)}$, (\ref{eqprtsb}) corresponds to 

\begin{equation}\label{eqreform}
 \PR_t \Big( \onpiv_{n - d(\phi,e^-)} \, \Big\vert \, \cross_{n - d(\phi,e^-)} \cap \seev_{n - d(\phi,e^-)}^c   \Big) \geq 
\PR_t \Big( \offpiv_{n - d(\phi,e^-)} \, \Big\vert \, \cross_{n - d(\phi,e^-)} \cap \seev_{n - d(\phi,e^-)}^c   \Big) \, .
\end{equation}

However, this is nothing other than Proposition~\ref{propineqone}.  
In this way, we see that the case of crossing with bottleneck reduces to that of crossing without bottleneck, and  Proposition~\ref{propineqtwo} indeed reduces to Proposition~\ref{propineqone}.

\medskip

\noindent{\bf Proof of Proposition~\ref{propineqtwo}.}
As we have noted, it is enough to verify~(\ref{eqprtsb})  for each $(e,h) \in E(\tree_{n-1}) \times [0,1)$. (Incidentally, 
note the $n$-dependence which is made explicit here. It is because we will reexpress the condition~(\ref{eqprtsb}) by another condition involving a different value of $n$ that we indicate the $n$-dependence of such quantities as $\cross$ and $\singlebar$ in this proof.)
\begin{definition}
For $v \in V(\tree)$, as a counterpart to the definition of descendent tree $\desctree{v}$, we define $\abovetree{v}$,     
\hfff{abovetree}
the ``tree $\tree$ above $v$''. This is the subtree of $\tree$ induced by all elements of $V(\tree)$ that are not strict descendents of $v$.

The set $V(\tree)$ may be labelled  by finite strings of symbols in $\{0,\ldots,\degg - 1 \}$. 
Concatenation of these labels
provides  an ordered addition operation on $V(\tree)$.
We extend the operation by setting $e +v = ( e^+ + v,e^- +v) \in E(\tree)$ for 
$e \in E(\tree)$ and
$v \in V(\tree)$. For $(v,h) \in V(\tree) \times [0,1)$ and for $b = (e',s) \in E(\tree_{[v]}) \times [0,1)$, define the $(v,h)$-shift $b^{(v,h)}$ of $b$ to be the bar $\big(e' - v,(s-h) \, {\rm mod} \, 1\big) \in E(\tree) \times [0,1)$. For a given bar set  $\ladder_0 \subseteq E(\tree_{[v]}) \times [0,1)$, define its 
$(v,h)$-shift $\ladder_0^{(v,h)}  \subseteq E(\tree) \times [0,1)$ to be $\big\{ b^{(v,h)} : b \in \ladder_0 \big\}$.
\end{definition}
\begin{definition}
For $m \in \N^+$, let $\condladderm$ denote the joint distribution of $\ladder \cap \big( E(\tree_m) \times [0,1) \big)$ and $\added_m \in  E(\tree_m) \times [0,1)$
under  $\PR_t\big( \cdot \big\vert \cross_m \cap \nse_m \big)$, where the law of $\added_m$ under $\PR_t$ is normalized Lebesgue measure on $E(\tree_m) \times [0,1)$ and $\cross_m$ and $\seev_m$ are the associated crossing and bottleneck events.
\end{definition}
\begin{lemma}\label{lemseev}
Let $(e,h) \in E(\tree_{n-1}) \times [0,1)$. For given $\ladder' \subseteq E(\tree) \times [0,1)$, write $\ladder'_{e^-} = \ladder' \cap \big( E(\tree_{[e^-]}) \times [0,1) \big)$. Under $\PR_t$, write $\ladder_{e^-}^{(e^-,h)}$ in place of $\big( \ladder_{e^-} \big)^{(e^-,h)}$. 
Then  the conditional joint law of $\ladder_{e^-}^{(e^-,h)}  \cap \big( E(\tree_{n - d(\phi,e^-)}) \times [0,1) \big)$
and $\added_n^{(e^-,h)}$ under $\PR_t \big( \cdot \big\vert \cross_n \cap \nonesc_n \cap \big\{ \singlebarn = (e,h) \big\}  \big)$ equals 
 $\PR^{\cross_{n - d(\phi,e^-)} \cap \seev_{n - d(\phi,e^-)}^c}_{t,\big(\ladder,\added_{n - d(\phi,e^-)}\big)}$.
\end{lemma}
The proof needs some definitions. 
%The first of these describes for a given choice of $\ladder$ the set of locations where the added bar may be placed in order that there be crossing without bottleneck.
\begin{definition}\label{defcrossed}
Let $\laddergiven \subset E(\tree) \times [0,1)$ be a given bar collection.
Let the set $\viloc_n(\laddergiven) \subseteq   E(\treen) \times [0,1)$ \hfff{viloc} of {\em viable} bar {\em locations} be such that $b \in \viloc_n(\laddergiven)$ if and only if both of the following conditions apply:
\begin{itemize}
 \item the meander $X^{\laddergiven}$ visits at least one joint of $b$ before time $\reach_n^{\laddergiven}$; 
 \item every edge in the path $P_{\phi,E(b)^+}$ supports at least two bars in $\laddergiven$.
\end{itemize}
Note that here $\reach_n^{\laddergiven}$ may or may not be finite.
\end{definition}
The key property of $\viloc_n$ is the following. Let $\laddergiven \subseteq E(\treen) \times [0,1)$. Conditionally on $\PR_t$ given $\ladder = \laddergiven$, 
\begin{equation}\label{eqkeyadded}
\added_n \in \viloc_n(\laddergiven) \; \, \textrm{if and only if} \; \, \cross_n \cap \seev_n^c \, .
\end{equation}
Here is the corresponding definition for a meander beginning from $(e^-,h)$ in place of $(\phi,0)$:
\begin{definition}\label{defviloces}
Let $(e,h) \in E(\tree_{n-1}) \times [0,1)$ and suppose given a bar collection $\laddergiven \subset 
%\big( 
%E(\treen) \cap 
E(\tree_{[e^-]}) 
%\big) 
\times [0,1)$.
Let $\viloc_n^{[e,h]}(\laddergiven) \subseteq \big(  E(\treen) \cap E(\tree_{[e^-]}) \big) \times [0,1)$ denote the set of $b \in \laddergiven$  such that both of the following conditions apply:
\begin{itemize}
 \item $X^{\laddergiven}_{(e^-,h)}$ visits at least one joint of $b$ before visiting $\scv_n \times [0,1)$; 
 \item every edge in the path $P_{e^-,E(b)^+}$ supports at least two bars in $\laddergiven$.
\end{itemize}
Here, the hitting time of $X^{\laddergiven}_{(e^-,h)}$ on  $\scv_n \times [0,1)$ may or may not be finite.
\end{definition}
\noindent{\bf Proof of Lemma \ref{lemseev}.}
Note that the following conditions are each necessary for $\cross_n \cap \nonesc \cap \big\{ \singlebar = (e,h) \big\}$:
\begin{itemize}
 \item
 $e$ supports exactly one bar in $\ladder$, this being $(e,h)$;
  \item $X^\ladder$ visits  $(e^+,h)$ before $\scv_n \times [0,1)$;
  \item $X^\ladder_{(e^+,h)}$ visits  $(\phi,0)$ before $\scv_n \times [0,1)$;
  \item $\added_n \in E(\tree_{[e^-]}) \times [0,1)$.
\end{itemize}
Note also that each of these events is measurable with respect to $(\ladder \cup \added_n) \cap \big(  E(\tree^{[e^-]}) \times [0,1) \big)$.
 Condition $\PR_t$ on their intersection, and 
note that the event $\cross_n \cap \nonescn \cap \big\{ \singlebarn = (e,h) \big\}$ is conditionally equal to $\added_n \in \viloc_n^{[e^-,h]}(\ladder_{e^-})$. 
Thus, under $\PR_t$ given $\cross_n \cap \nonesc_n \cap \{ \singlebarn = (e,h) \} \big)$, the conditional distribution of $\big(\ladder_{e^-},\added_n \big)$
is equal to an independent Poisson-$t$ random variable~$\ladder_{e^-}$ on $E(\tree_{[e^-]}) \times [0,1)$ and a Lebesgue-distributed element of  
$\big( E(\tree_{[e^-]}) \cap E(\tree_n) \big) \times [0,1)$ conditioned on $\added_n \in \viloc_n^{(e^-,h)}(\ladder_{e^-})$.
Note that the event $\added_n \in \viloc_n^{[e^-,h]}(\ladder_{e^-})$ coincides with 
$\added_n^{(e^-,h)}  \in \viloc_{n - d(\phi,e^-)}\big(\ladder_{e^-}^{(e^-,h)}\big)$.
That is,  under $\PR_t \big( \cdot \big\vert \cross_n \cap \nonesc_n \cap \{ \singlebarn = (e,h) \} \big)$, 
$\big( \ladder_{e^-}^{(e^-,h)},\added_n^{(e^-,h)} \big)$ is distributed as
$\big(\ladder,\added_{n - d(\phi,e^-)}\big)$ under $\PR_t$ given 
$\added_{n - d(\phi,e^-)}  \in \viloc_{n - d(\phi,e^-)}(\ladder)$; thus, the conditional law of  $\big(\ladder_{e^-}^{(e^-,h)},\added_n^{(e^-,h)}\big)$ is $\PR^{\cross_{n - d(\phi,e^-)} \cap \seev_{n - d(\phi,e^-)}^c}_{t,\big(\ladder,\added_{n - d(\phi,e^-)}\big)}$. \qed

\begin{lemma}\label{lemseevplus}
Let $(e,h) \in E(\tree_{n-1}) \times [0,1)$ be given. Let 
$\eta^\ladder \in [0,\infty]$ be the $\PR_t$-random variable given by $\inf \big\{ s > 0: X^\ladder_{(e^-,h)}(s) \in \{(e^-,h)\} \cup (\scv_n \times[0,1)) \big\}$; let $\eta^{\ladder \cup \added_n}$ denote the analogous stopping time for the process $X^{\ladder \cup \added_n}_{(e^-,h)}$. 
Then, 
under $\PR_t \big( \cdot \big\vert \cross_n \cap \nonesc_n \cap \big\{ \singlebarn = (e,h) \big\}  \big)$,
$\reach_n^\ladder = \infty$ if and only if $X_{(e^-,h)}^\ladder(\eta^\ladder) = (e^-,h)$, and 
$\reach_n^{\ladder \cup \added_n} = \infty$ if and only if $X_{(e^-,h)}^{\ladder \cup \added_n}(\eta^{\ladder \cup \added_n}) = (e^-,h)$.
\end{lemma}
\noindent{\bf Proof.} Under the law in question, $X^\ladder$ crosses $\singlebarn$ without having reached $\scv_n \times [0,1)$ (because this crossing must happen before $\hit^\ladder_{\added_n}$ which itself is before~$\reach_n^\ladder$).
After the crossing, $X^\ladder$ follows the trajectory of $X^\ladder_{(e^-,h)}$ for the duration $\eta^\ladder$, either ending up in $\scv_n \times [0,1)$ and thus realizing $\reach_n^\ladder < \infty$, or returning to $(e^-,h)$ and then 
pursuing the trajectory of $X^\ladder_{\singlebarn^+}$. In the latter case, $\nonesc_n$ ensures that $X^\ladder$ returns to $(\phi,0)$ before time $\reach_n^\ladder$, forcing this meander into a periodic trajectory and ensuring that $\reach_n^\ladder  = \infty$. Likewise for $X^{\ladder \cup \added_n}$. \qed

Lemmas \ref{lemseev} and \ref{lemseevplus} may be applied to reformulate (\ref{eqprtsb}) in the form~(\ref{eqreform}).
Since $n - d(\phi,e^-) \geq 1$, we may apply Proposition~\ref{propineqone} to find that, when its hypotheses hold, so does~(\ref{eqprtsb}).   This completes the proof of Proposition~\ref{propineqtwo}. \qed
\end{section}

\bibliographystyle{plain}

\bibliography{treecyclesbib}

\end{document}